\RequirePackage{fix-cm}
\documentclass[smallextended]{svjour3}       
\smartqed  
\usepackage{amsmath}
\usepackage{amsfonts}
\newtheorem{THEOREM}{Theorem}
\newtheorem{LEMMA}{Lemma}
\newtheorem{REMARK}{Remark}
\newtheorem{COROLLARY}{Corollary}
\begin{document}

\title{Geometry of Orlicz spaces equipped with norms generated by some lattice norms in $\mathbb{R}^{2}$
}


\author{Yunan Cui,     Henryk Hudzik,
        Haifeng Ma 
}


\institute{Yunan Cui \at
              Department of Mathematics, Harbin University of Science and Technology
Harbin,150080,China \\
              \email{cuiya@hrbust.edu.cn}           
           \and
           Henryk Hudzik \at
              Faculty of Economics and Information Technology, The State University of Applied Sciences in Plock, Nowe Trzepowo 55, 09-402 Plock, Poland
and Faculty of Mathematics and Faculty of Computer Science, Adam Mickiewicz University, Umultowska 87, 61-614 Poznan, Poland \\
              \email{hudzik@amu.edu.pl}
            \and
             Haifeng Ma \at
              School of Mathematical Science, Haribin Normal University, Harbin, 150025, China \\
              \email{haifengma@aliyun.com}
}

\date{}

\maketitle

\begin{abstract}
In Orlicz spaces generated by convex Orlicz functions a family of norms generated by some lattice norms in  $\mathbb{R}^{2}$ are defined and studied. This family of norms includes the family of the p-Amemiya norms ($1\leq p\leq\infty$) studied in [10-11], [14-15] and [20]. Criteria for strict monotonicity, lower and upper local uniform monotonicities and uniform monotonicities of Orlicz spaces and their subspaces of order continuous elements, equipped with these norms, are given in terms of the generating Orlicz functions, and the lattice norm in $\mathbb{R}^{2}$. The problems of strict convexity and of the existence of order almost isometric as well as of order isometric copies in these spaces are also discussed.
\end{abstract}

\section{Introduction}
\label{intro}
Let $p(\cdot)$ be a lattice norm in $\mathbb{R}^{2}$ such that $p((1,0))=1$. Let $\Phi:\mathbb{R}\rightarrow\mathbb{R}_{+}:=[0,+\infty)$ be an Orlicz function, that is, $\Phi$ is a nonzero function vanishing at zero, $\Phi$ is convex and even. Let us define $a(\Phi):=sup\{u\geq0:\Phi(u)=0\}$.\\
In the following $(\Omega,\Sigma,\mu)$ is a $\sigma$-finite complete measure space and $L^{0}=L^{0}( \Omega,\Sigma,\mu)$ is the space of all (equivalence classes of)  $\Sigma$-measurable functions $X:\Omega\rightarrow\mathbb{R} $ (where functions $x$ and $y$ equal $\mu$-a.e. in $\Omega$ belong to the same class of equivalence (we simply say that they are identified).\\
Given any Orlicz function $\Phi$ we define on the space $L^{0}( \Omega,\Sigma,\mu)$ the functional \\
\[I_{\Phi}(x)=\int\limits_{\Omega}\Phi(x(t))d\mu .\]
It is easy to see that the functional $I_{\Phi}$ has the following properties:\\
a) $I_{\Phi}:L^{0}( \Omega,\Sigma,\mu)\rightarrow\mathbb{R}^{e}_{+}:=\mathbb{R_+}\cup\{+\infty\}$,\\
b) $I_{\Phi}$ is convex,\\
c) $I_{\Phi}$ is even,\\
d) $I_{\Phi}(0)=0$ and if $x\in L^{0}( \Omega,\Sigma,\mu)\backslash \{0\}$, then $I_{\Phi}(\lambda x)\neq0$ for some $\lambda>0$.\\
and it is called the convex modular (see[39]). We are interested in introducing a norm generated by the functional $I_{\Phi}$ in the suitable biggest subspace of $L^{0}( \Omega,\Sigma,\mu)$. This subspace is called the Orlicz space, denoted by $L^{\Phi}=L^{\Phi}(\Omega,\Sigma,\mu)$ and defined by (see[5],[31-32],[35-37],[39-40])\\
\[ L^{\Phi}(\Omega,\Sigma,\mu):=\{x\in L^{0}( \Omega,\Sigma,\mu):I_{\Phi}(\lambda x)<\infty \,for \, some \, \lambda\in (0,+\infty)\}.\]
Let us demote by $A_{\Phi}(1)$  the modular unit ball, that is, \\
\[A_{\Phi}(1)=\{x\in L^{0}( \Omega,\Sigma,\mu):I_{\Phi}(x)\leq1\}.\]
Since the Orlicz function $\Phi$ is absolutely convex, that is,\\
\[\Phi(\alpha u+\beta \nu)\leq | \alpha| \Phi(u)+|\beta|\Phi(\nu) \, \]
for all u, $\nu\in \mathbb{R}$ and all $\alpha,\beta\in \mathbb{R}$ with $| \alpha|+|\beta|\leq1$, we obtain absolute convexity of the functional $I_{\Phi}$ and, in consequence, also the absolute convexity of the set $A_{\Phi}(1)$.
The Minkovski functional generated by the set $A_{\Phi}(1)$ can be defined for these elements from $ L^{0}( \Omega,\Sigma,\mu)$ which are absorbing by  $A_{\Phi}(1)$. It is easy to see that the biggest subspace of $ L^{0}( \Omega,\Sigma,\mu)$, the elements of which are absorbed by $A_{\Phi}(1)$, is just the Orlicz space $ L^{\Phi}(\mu) $. Namely,\\
\[ \exists \lambda > 0 ~s.t.:~ \frac{x}{\lambda}\in A_{\Phi}(1)) \Leftrightarrow(x\in  L^{\Phi}(\mu)).\]
The Minkovski functional of the set $ A_{\Phi}(1)$ is called in the literature the Luxemburg norm, it is denoted by $\|\cdot\|_{\Phi}$ and defined by the formula (see[5],[36],[37]and[39]) \\
\[\|x\|_{\Phi}=inf\{ \lambda>0: I_{\Phi}(\frac{x}{\lambda})\leq1\}\,(\forall x\in L^{\Phi}(\mu) ).\]
The following family of the norms, called p-Amemiya norms, was already defined and used in Orlicz spaces for $0\leq p\leq\infty$:\\
\[\|x\|_{\Phi,p}=\inf\limits_{k>0}\frac{1}{k}(1+( I_{\Phi}(kx))^{p})^{\frac{1}{p}}\,~(\forall x\in L^{\Phi}(\mu) ),\]
where $1\leq p\leq\infty ~$(see [10-11],[14-15],[20]). The norm $\|\cdot\|_{\Phi,1}$  is the Orlicz norm $\|\cdot\|^{0}_{\Phi}$ which was defined by Orlicz in [41] by the formula\\
\[\|x\|^{0}_{\Phi}=sup\{ |\int\limits_{\Omega}x(t)y(t)d\mu|:y\in L^{0}( \Omega,\Sigma,\mu)~and~I_{\Phi^*}(y)\leq1  \}, \]
where $\Phi^{*}$ is the function complementary to $\Phi$ in the sense of Young, that is,\\
\[\Phi^{*}(u)=\sup\limits_{\nu\geq0}\{ |u|\nu-\Phi(\nu) \}.\]
For $p=\infty$, we have (see [24])\\
\[\|x\|_{\Phi,\infty}=   \lim\limits_{p\rightarrow\infty}\frac{1}{k}(1+( I_{\Phi}(kx))^{p})^{\frac{1}{p}}=\inf\limits_{k>0}\frac{1}{k}max(1,( I_{\Phi}(kx)))=\|x\|_{\Phi}.\]
It is well known that all the norms from the family $\{\|\cdot\|_{\Phi,p}\}_{p\in [1,\infty]}$ are equivalent and that the Luxemburg norm
$\|\cdot\|_{\Phi}= \|x\|_{\Phi,\infty}  $ is the smallest norm and the Orlicz norm $\|\cdot\|^{0}_{\Phi}= \|x\|_{\Phi,1}  $ is the biggest one. The Orlicz space $ L^{\Phi}(\mu)$ equipped with every norm from this family of norms is a Banach space, which is even the Banach function lattice, called also the K\H{o}the  space (see[31-32],[35]and[41]), which means that for any $p\in [1,+\infty]$, the space ($ L^{\Phi}(\mu)$, $\|\cdot\|_{\Phi,p} $) has the following properties:\\
$1^{\circ}$ For any $x\in L^{\circ}(\mu)$, $y\in L^{\Phi}(\mu)$, if $|x(t)|\leq|y(t)|$ for $\mu$ a.e. $t\in \Omega$, then $x\in L^{\Phi}$ and $\|x\|_{\Phi,p}\leq \|y\|_{\Phi,p}$,\\
$2^{\circ}$ There exist a function $x\in L^{\Phi}(\mu)$ such that $x(t)>0$ for any $t\in \Omega$.\\
The same properties has the space $(E^{\Phi}(\mu),\|\cdot\|_{\Phi,p} ) $ defined below. Let us recall that an element $x$ of a K\H{o}the space $(E,\,\|\cdot\|_{E})$ is said to be order continuous if for any sequence $\{x_{n}\}^{\infty}_{n=1}$ in $E$ such that $0\leq x_{n}(t)\leq |x(t)| $ for all   $n\in \mathbb{N}$ and $\mu$ a.e.$t\in T$, the condition $x_{n}(t)\rightarrow 0$ as $n\rightarrow\infty$ for $\mu$ a.e. $t\in \Omega$ implies that $\|x_{n}\|_{E}\rightarrow0$ as $n\rightarrow\infty$. The set of all order continuous elements in $E$ is denoted by $E_{a}$, and the space $(E_{a},\,\|\cdot\|_{E}) $ is again a K\H{o}the space. It is obvious that equivalent norms keep the order continuity property. It is well known that (see[35]and[41])\\
\[(L^{\Phi}(\mu))_{a}=E^{\Phi}(\mu),\]
where \\
\[E^{\Phi}(\mu):=\{x\in L^{0}( \Omega,\Sigma,\mu):I_{\Phi}(\lambda x)<\infty \, for\, any \, \lambda>0 \}.\]
In this paper we will introduce a new family of norms in the Orlicz space  $L^{\Phi}(\mu)$. Namely, given any lattice norm $p(\cdot)$ in $\mathbb{R}^{2}$ such that $p((1,0))=1$, we define the following functional in $L^{\Phi}(\mu)$:\\
\[\|x\|_{\Phi,p(\cdot)}:=\inf\limits_{k>0}\frac{1}{k}p((1, I_{\Phi}(kx)))\,~(\forall x\in L^{\Phi}(\mu) ).\]
We will prove that such functionals are norms in $L^{\Phi}(\mu)$. Of course, these norms are equivalent each others. We will work on criteria for strict convexity and various their monotonicity properties (strict monotonicity, lower and upper local uniform monotonicity and uniform monotonicity) as well as on order almost isometric copies of $l^{\infty} $ and order isometric copies of $l^{\infty} $.\\
We need to define all others notions that will be used in this paper. A Banach lattice $X=(X,\leq,\|\cdot\|)$, for the definition of which we refer to [2], [31], [35] and [41] is said to be strictly monotone if for any $x,y\in X$ such that $0\leq x\leq y$ and $x\neq y$ we have $\|x\|<\|y\|$. By the homogenity of the norm $\|\cdot\|$, we can restrict ourselves in this definition to $y\geq0$ satisfying $\|y\|=1$. Let us denote by $X_{+}$ the positive cone in $X$ that is the set of all $x\in X$ such that $x\geq 0$. In our definitions below $X$ always denotes a Banach lattice $(X,\leq,\|\cdot\|)$.
$X$ is said to be uniformly monotone (see [2] and [33]) if for any $ \varepsilon \in (0,1)$ there exists $\delta(\varepsilon)\in (0,1)$ such that if $x,y\in X$, $0\leq x\leq y$; $\|x\|\geq\varepsilon $ and $\|y\|=1$, then $\|y-x\|\leq1-\delta(\varepsilon)$. The biggest function $\delta_X:(0,1)\rightarrow(0,1)$ with this property, that is, the function
\[ \delta_{X}(\varepsilon)=inf\{1-\|y-x\|: 0\leq x\leq y;\|x\|\geq\varepsilon,\|y\|=1 \} \]
is called the modulus of monotonicity of $X$ (see [2]) and for the properties of $\delta_{X}(\cdot)$ also ([23]). It is known (see [33]) that $X$ is uniformly monotone if and only if for any $\varepsilon>0$ there exists $\sigma(\varepsilon)>0$ such that for any $x,y\in X_{+}$ such that $\|x\|\geq\varepsilon $ and $\|y\|=1$ there holds $\|y+x\|\geq1+ \sigma(\varepsilon)$. $X$ is said to be lower (upper) locally uniformly monotone if for any $y\in X_{+}$ with $\|y\|=1$ and any $ \varepsilon\in (0,1)$ (resp. any $\varepsilon>0$) there exists $\delta(y,\varepsilon)\in (0,1)$(resp. $\sigma(y,\varepsilon)>0$ ) such that for any $x\in X$ satisfying $0\leq x\leq y$ and $\|x\|\geq\varepsilon$ (resp. $x\geq 0$ with $\|x\|\geq\varepsilon$), we have $\|y-x\|\leq1- \delta(y,\varepsilon)$ (resp. $\|y+x\|\geq1+ \sigma(y,\varepsilon)$). For the definition of these two properties see [2,4,25].\par
It is obvious that $\Phi$ vanishes only only at 0 iff $a(\Phi)=0$. For any Orlicz function $\Phi$ we say that it satisfies condition $\Delta_{2}(\mathbb{R}_{+})$ ($\Phi\in \Delta_{2}(\mathbb{R}_{+})$ for short) if there exists $K>0$ such that $\Phi(2u)\leq K\Phi(u)$ for any $u\geq0$. We say that $\Phi$ satisfies condition $\Delta_{2}$ at infinity ($\Phi\in \Delta_{2}(\infty)$ for short)
if there are positive constants $u_{0},K$ such that $\Phi(2u)\leq K\Phi(u)$ for all $u\geq{u_{0}}$. We say that $\Phi$ satisfies condition $\Delta_{2}$ at zero ($\Phi\in \Delta_{2}(0)$ for short) if there exist two positive constants $u_{0},K$ such that $ \Phi(2u)\leq K\Phi(u)$ for all $u\in [0,u_{0}]$. It is easy to see that $\Delta_{2}(\mathbb{R}_{+})$ if and only if $\Phi\in \Delta_{2}(\infty)$ and $\Phi\in \Delta_{2}(0)$.\par
The $\Delta_{2}$-condition for $\Phi$ should be defined suitably to the measure space $(\Omega,\Sigma,\mu)$ in such a way that the corresponding Orlicz space $(L^{\Phi}(\mu),\|\cdot\|_{\Phi} )$ is order continuous. We know that suitable $\Delta_{2}$-condition for the couple $(\Phi,(\Omega,\Sigma,\mu))$ is the following:\\
a) condition $\Delta_{2}(\mathbb{R}_{+})$ for $\Phi$ if $(\Omega,\Sigma,\mu)$ is infinite and non-atomic.\\
b) condition $\Delta_{2}(\infty)$ for $\Phi$ if $(\Omega,\Sigma,\mu)$ is finite and non-atomic.\\
c) condition $\Delta_{2}(0)$ for $\Phi$ if $\Omega=\mathbb{N}$, $\Sigma=2^{\mathbb{N}}$ and $\mu$ is the counting measure on $2^{\mathbb{N}}$.\\
It is obvious that any $\sigma$-finite measure space $ (\Omega,\Sigma,\mu)$ can be represented as the direct sum of two measure spaces
$(\Omega_{n-a},\Sigma\cap\Omega_{n-a},\mu|_{\Sigma\cap\Omega_{n-a}})\oplus(\Omega_{a},\Sigma\cap\Omega_{a},\mu|_{\Sigma\cap\Omega_{a}})$, where $\Omega_{a}$ is the set of all atoms for $\mu$ in $\Sigma$ and $ \Omega_{n-a}=\Omega\backslash \Omega_{a}$.\par
If $\Omega_{a}$ is finite and $\mu(\Omega_{n-a})>0$, then the suitable condition $\Delta_{2}$ for $\Phi$ is the $\Delta_{2}$-condition for the non-atomic measure space defined above. If $\mu(\Omega_{n-a})=0$ and $\Omega_{a}=\mathbb{N}$, $\Sigma=2^{\mathbb{N}}$ and $\mu$ is the counting measure on $2^{\mathbb{N}}$, then the suitable condition $\Delta_{2}$ is the condition $\Delta_{2}(0)$. If $\mu(\Omega_{n-a})>0$, $\Omega_{a}=\mathbb{N}$, $\Sigma=2^{\mathbb{N}}$ and $\mu$ is the counting measure on $2^{\mathbb{N}}$, then the suitable $\Delta_{2}$-condition for $\Phi$ is the conjunction of the suitable $\Delta_{2}$-condition for the non-atomic measure space and of the condition $\bigtriangleup_2(0)$. In our paper, we always assume that all atoms have the measure 1 and we identify the atoms with the singletons $\{n\}$, where $n\in \mathbb{N}$ (the set of all natural numbers).\par
Monotonicity properties of Banach lattices have applications in the dominated best approximation (see[33],[25],[8],[17] and [14]) and in the fixed point theory (see[14],[17] and [20]). They are also strongly related to the complex  rotundity  properties (see [30]). For these reasons monotonicity  properties were investigated in various classes of function spaces. Namely, in [13], [25-27], [33], [34] and [16] for Musielak-Orlicz spaces, in [21] for Lorentz spaces, in [19] for Orlicz-Lorentz spaces, in [4], [14], [16], [23], [38] for Orlicz spaces, in [29-30] for Calder$\grave{o}$n-Lozanovski$\check{\i}$ spaces, in [12] for Ces$\grave{a}$ro-Orlicz sequence spaces, for Orlicz-Sobolev space in [8].
Relationships between monotonicity properties and rofundity properties as well as between monotonicity properties and orthogonal monotonicity properties in  K\H{o}the    spaces were studied in [22]. In abstract Banach lattices relationships between monotonicity properties and dominated best approximation problems we studied in [6], [14] and [17].\par

Problems on estimates or calculations of the characteristic of monotonicity in Orlicz spaces and Orlicz-Lorentz spaces were studied in [16], [19] and [23]. Applications of the monotonicity properties and the ergodic theory in Banach lattices were studied in [1].\par
\begin{THEOREM}
For any lattice norm $p(\cdot)$ on $\mathbb{R}^{2}$ such that $p((1,0))=p((0,1))=1$ we have the inequality\\
\[\forall (u,\nu)\in \mathbb{R}^{2}:\quad max(|u|,|\nu|)\leq p((u,\nu))\leq |u|+|\nu|\]\\
that is the smallest (resp.the biggest) lattice norm $p(\cdot)$ among these ones with $p((1,0))=p((0,1))$ is the $l^{\infty}-$norm (resp. $l^{1}-$norm).
\end{THEOREM}
{\it Proof.}
Let us take any $(u,\nu)\in \mathbb{R}^{2}$. Then \\
\begin{eqnarray*}
p((u,\nu)) &=& p((|u|,|\nu|))\geq p((|u|,0))=p(|u|(1,0)) \\
&=& |u|p((1,0))=|u|
\end{eqnarray*}
and \\
\begin{eqnarray*}
p((u,\nu)) &=& p((|u|,|\nu|))\geq p((0,|\nu| ))=:p(|\nu|(0,1)) \\
&=& |\nu|p((0,1))=|\nu|,
\end{eqnarray*}
whence\\
\[p((u,\nu)) \geq max(|u|,|\nu|)=p_{\infty}((u,\nu)) .\]\\
On the other hand for any $(u,\nu)\in \mathbb{R}^{2}$, we have \\
\begin{eqnarray*}
p((u,\nu)) &=& p((|u|,|\nu|)) =p((|u|,0)+(0,|\nu| )) \\
&\leq& p((|u|,0))+p((0,|\nu| ))=|u|p((1,0))+|\nu|p((0,1 ))\\
&=& |u|+|\nu|=:p_1((u,\nu)).
\end{eqnarray*}
It is known that $\|\cdot\|_{\Phi,p_{\infty}(\cdot)}$ is equal to the Luxembourg norm and $\|\cdot\|_{\Phi,p_{1}(\cdot)}$  is equal to the Orlicz norm as well as to the Amemiya norm (see [24]). Geometry  of Orlicz spaces equipped with the $p-$Amemiya norm, that is, the norm
$ \|\cdot\|_{{p}(\cdot)} $, where $p((u,\nu)):=(|u|^{p}+|\nu|^{p})^{\frac{1}{p}}$ for any $(u,\nu)\in \mathbb{R}^{2}$, was considered in the papers [10-11], [14-15], the fixed point property in these spaces was studied in [14] and [20]. The dominated best approximation in these spaces was studied in [14].

\begin{THEOREM}
For any Orlicz function $\Phi$ and any lattice norm $p(\cdot)$ in the Eucklidean space $\mathbb{R}^{2}$ the functional $\|\cdot\|_{\Phi,p(\cdot)}$ is a norm in $L^{\Phi}(\mu)$.
\end{THEOREM}
{\it Proof.} We have $I_{\Phi}(k0)=0$ for any $k>0$, so \\

\[{\|0\|_{\Phi,p(\cdot)}=\inf\limits_{k>0}\frac{1}{k}p((1,0))}=\inf\limits_{k>0}\frac{1}{k}=0. \]
Let us assume that $x\in L^{\Phi}(\mu)\backslash\{0\}$. Since $\Phi$ is a nonzero function, that is, there exists $u_{0}>0$ such that $\Phi(u_{0})>0$, so
there exists $k_{0}>0$ such that $I_{\Phi}(k_{0}x)>1$. Then
\begin{eqnarray*}
\|x\|_{\Phi,p(\cdot)} &=& \min(\inf\limits_{0<k\leq k_{0}}\frac{1}{k}p((1, I_{\Phi}(kx))),\inf\limits_{k\geq k_{0}}\frac{1}{k}p((1,I_{\Phi}(kx))))\\
 &\geq& \min(\inf\limits_{0<k\leq k_{0}}\frac{1}{k}p((1, 0)),\inf\limits_{k\geq k_{0}}\frac{1}{k}p((0,I_{\Phi}(kx))))\\
 &\geq& \min(\frac{1}{k_{0}},\inf\limits_{k\geq k_{0}}p((0,\frac{1}{k}I_{\Phi}(kx))))\\
 &=& \min(\frac{1}{k_{0}},p((0,\frac{1}{k_{0}}I_{\Phi}(k_{0}))))\\
 &>& 0.
\end{eqnarray*}
Now, we will show that the functional $\| \cdot\|_{\Phi,p(\cdot)}$ is absolutely homogeneous. Let us take any $x\in L^{\Phi}(\mu)$ and any $\lambda>0$. If $x=0$, then $\lambda x=0$ for any $\lambda\in \mathbb{R}$, whence

$$\|\lambda 0\|_{\Phi,p(\cdot)}=\|0\|_{\Phi,p(\cdot)}=0=\lambda\cdot0=\lambda\|0\|_{\Phi,p(\cdot)}.$$
So let us assume that $x\neq0$. Then

\begin{eqnarray*}
\|\lambda x\|_{\Phi,p(\cdot)} &=& \inf\limits_{k>0}\frac{1}{k}p((1, I_{\Phi}(k |\lambda| x)))\\
 &=& \inf\limits_{k>0}\frac{|\lambda|}{|\lambda|k}p((1, I_{\Phi}(k \lambda x)))\\
  &=& |\lambda|  \inf\limits_{k>0}\frac{1}{|\lambda|k}p((1, I_{\Phi}(k |\lambda| x)))\\
 &=&  |\lambda| \| x\|_{\Phi,p(\cdot)}
\end{eqnarray*}
Finally, we will show that the functional $\|\cdot\|_{\Phi,p(\cdot)}$ satisfies the triangle inequality. Let us take arbitrarily $x,y\in L^{\Phi}(\mu)$. If at least one element among $x$ and $y$ is equal to zero function, then the triangle inequality is obvious. So assume that $x\neq 0$ and $y\neq 0$. Let us take any $\varepsilon>0$. There exists constants $\lambda>0$ and $l>0$ such that \\
\[  \frac{1}{\lambda}p((1, I_{\Phi}(\lambda x)))\leq  \|x\|_{\Phi,p(\cdot)}+\varepsilon,\]
\[  \frac{1}{l}p((1, I_{\Phi}(l y)))\leq  \|y\|_{\Phi,p(\cdot)}+\varepsilon.\]
Then
\begin{eqnarray*}
\|x+y\|_{\Phi,p(\cdot)} &\leq& \frac{\lambda+l}{\lambda l}p((1, I_{\Phi}(\frac{\lambda l}{\lambda +l} (x+y))))\\
 &=& \frac{\lambda+l}{\lambda l}p((1, I_{\Phi}(\frac{ l}{\lambda +l}( \lambda x)+ \frac{ \lambda}{\lambda +l}( l y))))\\
 &=& \frac{\lambda+l}{\lambda l}p(( \frac{ l}{\lambda +l}+\frac{ \lambda}{\lambda +l}     , I_{\Phi}(\frac{ l}{\lambda +l}( \lambda x)+ \frac{ \lambda}{\lambda +l}( l y))))\\
 &\leq& \frac{\lambda+l}{\lambda l}p(( \frac{ l}{\lambda +l}+\frac{ \lambda}{\lambda +l},   \frac{ l}{\lambda +l} I_{\Phi}( \lambda x)+ \frac{ \lambda}{\lambda +l} I_{\Phi}( l y)))\\
 &=& \frac{\lambda+l}{\lambda l}p(( \frac{ l}{\lambda +l},   \frac{ l}{\lambda +l} I_{\Phi}( \lambda x)) +(\frac{ \lambda}{\lambda +l},\frac{ \lambda}{\lambda +l} I_{\Phi}( l y)))\\
 &\leq& \frac{\lambda+l}{\lambda l} \{\frac{ l}{\lambda +l}p(( 1,I_{\Phi}( \lambda x))) +\frac{ \lambda}{\lambda +l}p((1, I_{\Phi}( l y)))\}\\
 &=& \frac{1}{\lambda } p(( 1,I_{\Phi}( \lambda x))) +\frac{ 1}{l}p((1, I_{\Phi}( l y)))\\
&\leq&  \|x\|_{\Phi,p(\cdot)}+\|y\|_{\Phi,p(\cdot)}+2\varepsilon.
 \end{eqnarray*}
 By the arbitrariness of $\varepsilon>0$, we obtain the inequality
 \[ \|x+y\|_{\Phi,p(\cdot)} \leq  \|x\|_{\Phi,p(\cdot)}+\|y\|_{\Phi,p(\cdot)} , \]
 which finishes the proof of the theorem.
  \begin{LEMMA}
 If $p(\cdot)$ is a lattice norm in $\mathbb{R}^{2}$ and $\Phi$ be an Orlicz function satisfying the condition $\lim\limits_{u\rightarrow+\infty}(\Phi(u)/u)=+\infty$, then for any $x\in L^{\Phi}(\mu)\backslash\{0\}$ there exists $l\in (0,+\infty)$ such that
 \[\|x\|_{\Phi,p(\cdot)}=\frac{1}{l}p((1, I_{\Phi}(l x))).\]
 \end{LEMMA}
 {\it Proof.} Since $p(\cdot)$ is a lattice norm in $\mathbb{R}^{2}$, we have for any $u>0$ that
 \[ p((1,u))\geq p((0,u))=up((0, 1))\rightarrow\infty \quad as \quad u\rightarrow\infty .\]
 Since $x\neq0$, the condition $(\Phi(u)/u)\rightarrow+\infty$ as $u\rightarrow\infty$ implies that \\
 (1) \quad $\frac{1}{k}I_{\Phi}(kx)\rightarrow\infty$ as $k\rightarrow\infty$.\\
 Since $x\in L^{\Phi}(\mu)$, there exists $\lambda>0 $ such that $0<I_{\Phi}(\lambda x)<\infty$. Hence, and from condition (1) as well as by the fact that
 \[ \lim\limits_{k\rightarrow0^{+}} \frac{1}{k}p((1, I_{\Phi}(k x))) \geq \lim\limits_{k\rightarrow0^{+}} \frac{1}{k}p((1,0))=\lim\limits_{k\rightarrow0^{+}} \frac{1}{k}=+\infty, \]
 there exist positive constants $k_{0}$ and $k_{1}$ such that $k_{0}<k_{1}<+\infty$, $I_{\Phi}(k_{1}x)<\infty $, and
\[\|x\|_{\Phi,p(\cdot)}=\inf\limits_{ k_{0}\leq{k}\leq{k_{1}}}\frac{1}{k}p((1, I_{\Phi}(k x))).\]
The function $ f(k):=I_{\Phi}(k x))$ is convex and it has finite values on the compact interval $[k_{0},k_{1}]$, so it is continuous on this interval. In consequence, by continuity of the norm $p(\cdot)$, the function $g:[k_{0},k_{1}] \rightarrow\mathbb{R}_{+}$ defined by
\[g(k)=p((\frac{1}{k}, \frac{1}{k} I_{\Phi}(k x)))\]
is also continuous. Therefore, the desired number $l\in (0,+\infty)$ exists.
\begin{LEMMA}
For any Orlicz function $\Phi$, any lattice norm on $\mathbb{R}^{2}$ such that  $p((1,0))=p((0,1))=1$, if $x\in L^{\Phi}(\mu)$ is such that $I_{\Phi}(\lambda{x})=+\infty$ for any $\lambda>1$, then
\[1\leq\|x\|_{\Phi,p(\cdot)}\leq{1+I_{\Phi}(x)} \]
\end{LEMMA}
{\it Proof.} We have
 \begin{eqnarray*}
\| x\|_{\Phi,p(\cdot)} &=& \inf\limits_{k>0}\frac{1}{k}p((1, I_{\Phi}(k x)))\\
 &=& \min(\inf\limits_{0<k\leq1}\frac{1}{k}p((1,I_{\Phi}(kx))),\inf\limits_{k\geq 1}p((1,I_{\Phi}(kx))))\\
 &\geq& \min(\inf\limits_{0<k\leq1}\frac{1}{k}p((1,0)),p((1,I_{\Phi}(kx))))\\
 &=& \min(1,p((1,I_{\Phi}(kx))))\\
 &=& 1.
\end{eqnarray*}
 On the other hand
 \begin{eqnarray*}
\| x\|_{\Phi,p(\cdot)} &=& \inf\limits_{k>0}\frac{1}{k}p((1, I_{\Phi}(k x)))\\
 &\leq& p((1,I_{\Phi}(x)))=p((1,0)+( 0, I_{\Phi}(x)) )\\
 &\leq& p((1,0))+p(( 0, I_{\Phi}(x)) )\\
 &=& 1+ I_{\Phi}(x)p((0,1))=1+ I_{\Phi}(x).
\end{eqnarray*}

\begin{THEOREM}
Let $p(\cdot)$ be a norm on $\mathbb{R}^{2}$ as in Lemma 2. If $\Phi$ is an Orlicz function $\Phi$ which does not satisfy suitable $\Delta_{2}$-condition, then $(L^{\Phi}(\mu), \| \cdot\|_{\Phi,p(\cdot)})$ contains an order linearly almost isometric copy of $l^{\infty}$, that is, for any $\varepsilon>0$ there exists a linear nonnegative operator $P_{\varepsilon}:l^{\infty}\rightarrow L^{\Phi}(\mu)$ such that
\[ \|z\|_{\infty}\leq\|Pz\|_{\Phi,p(\cdot)}\leq(1+\varepsilon) \|z\|_{\infty}\quad (  \forall z\in l^{\infty} ).\]

\end{THEOREM}
{\it Proof.}
Under the assumptions on $\Phi$, given any $\varepsilon>0$, there exists a sequence $\{x_{n}\}^{\infty}_{n=1}$ in $L^{\Phi}(\mu)$ with pairwise disjoint supports and such that $I_{\Phi}(x_{n})\leq \varepsilon/2^{n}$, $x_{n}\geq0$ and  $I_{\Phi}( \lambda x_{n})=\infty$ for any $n\in \mathbb{N}$ and $\lambda>1$. Let us define the operator $P_{\varepsilon}$ on $L^{\Phi}(\mu)$ by the formula
\[P_{\varepsilon}z=\sum_{n=1}^{\infty}z_{n}x_{n}\quad (\forall z=\{z_{n}\}_{n=1}\in l^{\infty}),  \]
where the series is defined pointwisely for $t\in\Omega$. It is obvious that $P_{\varepsilon}$ is linear and nonnegative. There is no problem with the pointwise convergence of the series by pairwise disjointness of the supports of the element $x_{n}\in L^{\Phi}(\mu),$ whence for any $t\in\Omega$ there exists at least one $n\in N$ such that $t\in suup x_{n}$.\par
Let us note first that for any $k>0$ and $z\in l^{\infty}$, we have
\begin{eqnarray*}
p((1,I_{\Phi}(k\frac{P_{\varepsilon}z}{\|z\|_{\infty}})))&=&p((1,I_{\Phi}(k\frac{P_{\varepsilon}|z|}{\|z\|_{\infty}})))\\
&=& p((1,I_{\Phi}(k\frac{\sum_{n=1}^{\infty}|z_{n}|x_{n}}{\|z\|_{\infty}})))\leq p((1,I_{\Phi}(k\sum_{n=1}^{\infty}x_{n})))\\
&=& p((1,\sum_{n=1}^{\infty}I_{\Phi}(kx_{n}))).
\end{eqnarray*}
In consequence, by Lemma 2, we have for any $z=\{z_{n}\}_{n=1}^{\infty}\in l^{\infty}$,
\begin{eqnarray*}
\|\frac{P_{\varepsilon}z}{\|z\|_{\infty}}\|_{\Phi,p(\cdot)}&=&\inf\limits_{k>0} \frac{1}{k}p((1,I_{\Phi}(k\frac{P_{\varepsilon}z}{\|z\|_{\infty}})))\\
&\leq&  p((1,I_{\Phi}(\frac{P_{\varepsilon}z}{\|z\|_{\infty}})))\\
&\leq& p((1,I_{\Phi}(\sum_{n=1}^{\infty}x_{n})))\\
&\leq& 1+\varepsilon,
\end{eqnarray*}
whence $\|P_{\varepsilon}z\|_{\Phi,p(\cdot)}\leq(1+\varepsilon)\|z\|_{\infty }$ for any $z\in l^{\infty}$.\par
On the other hand, since for any $\lambda>0$ there exist $n_{\lambda}\in N$ such that $\lambda|z_{n_{l}}|/\|z\|_{\infty}>1$,we have

\[I_{\Phi}(\frac{P_{\varepsilon}\lambda z}{\|z\|_{\infty}})\geq I_{\Phi}(\frac{\lambda|z_{n_{\lambda}}|}{\|z\|_{\infty}}x_{n_{\lambda}})=\infty,\]\\
whence, by Lemma 2,\\
 \[\|\frac{P_{\varepsilon}\lambda z}{\|z\|_{\infty}}\|_{\Phi,p(\cdot)}=\|\frac{\lambda P_{\varepsilon}z}{\|z\|_{\infty}}\|_{\Phi,p(\cdot)}\geq1,\]\\
that is,
\[\|P_{\varepsilon}z\|_{\Phi,p(\cdot)}\geq \frac{\|z\|_{\infty}}{\lambda}.\]
By the arbitrariness of $\lambda>1$, we have $\|P_{\varepsilon}z\|_{\Phi,p(\cdot)}\geq \|z\|_{\infty}$ for any $z\in l^{\infty}$, which finishes the proof.\\
\begin{THEOREM}
Let $p(\cdot)$ be a norm in $\mathbb{R}^{2}$ such as in Lemma 2 and $\Phi$ be an Orlicz function with $a(\Phi):=sup\{u\geq0:\Phi(u)\}>0$. Then in both cases, a non-atomic infinite measure space as well as the case of the counting measure on $2^{N}$, the Orlicz space $(L^{\Phi}(\mu),\|\cdot\|_{\Phi,p(\cdot)})$
contains a linearly order isometric copy of $l^{\infty}$.\\
\end{THEOREM}
{\it Proof.} Under the assumptions on the measure space, there exists a sequence $\{A_{n}\}_{n=1}^{\infty}$ of pairwise disjoint set with $\mu(A_{n})=+\infty$
 for any $n\in N$. Let us define \\
 \[x_{n}=a(\Phi)\mathcal{X}_{A_{n}}, ({\forall}n\in N), x:=\sum\limits_{n=1}^{\infty}x_{n}=\sup\limits_{n\in N} x_n,\] \\
 where the series is defined pointwisely (no problem with it's pointwise convergence because of pairwise disjointness of the sets $A_{n}$). It is obvious that $I_{\Phi}(x)=0$ and $I_{\Phi}(x_{n})=0$ as well as that $I_{\Phi}(\lambda x)=I_{\Phi}(\lambda x_{n})=+\infty$ for any $n\in N$ and $\lambda>1$. Moreover,\\
 \[\begin{aligned}
\parallel x_{n}\parallel_{\Phi,p(\cdot)}&=min(\inf_{0<k\leq1}\frac{1}{k}p((1,I_{\Phi}(kx)))), \inf_{k\geq1}\frac{1}{k}p((1,I_{\Phi}(kx_{n})))\\
&=min(\inf_{0<k\leq1}\frac{1}{k}, p((1,I_{\Phi}(x_{n}))))\\
&=min(1, 1)=1 ~({\forall}n\in N).
\end{aligned}\]
In the same way, we can prove that $\parallel x\parallel_{\Phi,p(\cdot)}=1$. Let us define the following operator on $l^{\infty}$:\\
\[Pz=\sum\limits_{n=1}^{\infty}z_{n}x_{n}\,\,\,\,     ({\forall}z=\{z_{n}\}\in l^{\infty}).
\]\\
Let us first note that $P: l^{\infty}\longrightarrow L^{\Phi}(\mu).$ Namely, $I_{\Phi}(\frac{Pz}{\parallel z\parallel_{\infty}})\leq I_{\Phi}(\sum\limits_{n=1}^{\infty}x_{n})=I_{\Phi}(x)=0$, whence $Pz\in L^{\Phi}(\mu)$ for any $z\in l^{\infty}$. Moreover,\\
\[\begin{aligned}
\parallel\frac{Pz}{\parallel z\parallel_{\infty}}\parallel_{\Phi,p(\cdot)}&=\inf_{k>0}\frac{1}{k}p((1, I_{\Phi}(k\frac{Pz}{\parallel z\parallel_{\infty}})))\\
&\leq \inf_{k>0}\frac{1}{k}p((1, I_{\Phi}(kx)))\\
&=\parallel x\parallel_{\Phi,p(\cdot)}\\
&=1,
\end{aligned}\]\\
whence $\parallel Pz\parallel_{\Phi,p(\cdot)}\leq \parallel z\parallel_{\infty}$. On the other hand, given any $\lambda>1$, one can find $n_{\lambda}\in N$ such that $\lambda\mid z_{n_{\lambda}}\mid>\parallel z\parallel_{\infty}$, consequently,\\
\[\parallel\frac{\lambda Pz}{\parallel z\parallel_{\infty}}\parallel_{\Phi,p(\cdot)}=\parallel\frac{P(\lambda z)}{\parallel z\parallel_{\infty}}\parallel_{\Phi,p(\cdot)}>\parallel x_{n_{\lambda}}\parallel_{\Phi,p(\cdot)}=1,\]
whence $\parallel Pz\parallel_{\Phi,p(\cdot)}>\frac{\parallel z\parallel_{\infty}}{\lambda}$. By the arbitrariness of $\lambda>1$, we obtain that $\parallel Pz\parallel_{\Phi,p(\cdot)}\geq\parallel z\parallel_{\infty}$, which together with the opposite inequality proved already gives the equality $\parallel Pz\parallel_{\Phi,p(\cdot)}= \parallel z\parallel_{\infty}$ for any $z\in l^{\infty}$, which means that $P$ is an isometry. It is obvious that the operator $P$ is linear. Since the functions $x_{n}$ are non-negative, so $P$ is also non-negative, that is, $Pz\geq0$ for any $z\in l^{\infty}$, $z\geq0$. In consequence, the operator $P$ is a linear order isometry, which finishes the proof. \\
\begin{THEOREM} Let $p(\cdot)$ be a lattice norm in $\mathbb{R}^{2}$ which is strictly increasing on the vertical half-line $\{(1,u): u\in R_+\}$ in $\mathbb{R}^{2}$ and let $\Phi$ be a strictly convex Orlicz function. Let the couple $(\Phi, p(\cdot))$ satisfy the condition\\
(2) $({\forall}x\in L^{\Phi}(\mu)\backslash\{0\})$ $(\exists l\in (0, +\infty))$ $\parallel x\parallel_{\Phi, p(\cdot)}=\frac{1}{l}p((1, I_{\Phi}(lx))).$\\
Then the space $(L^{\Phi}(\mu), \parallel \cdot\parallel_{\Phi, p(\cdot)})$ is strictly convex.\\
\end{THEOREM}
{\it Proof.} Assume that $x,y\in S((L_{+}^{\Phi}(\mu),\parallel \cdot\parallel_{\Phi, p(\cdot)}))$ and $x\neq y$. Then $\Phi\circ kx\neq\Phi\circ ky$ for any $k\in(0, +\infty)$ because strict convexity of $\Phi$ implies that $\Phi$ is a $1-1$ function on $R_+$. We also have that for any $k\in(0, +\infty)$. Let $\lambda$, $l\in(0, +\infty)$ be such that\\
\[\parallel x\parallel_{\Phi,p(\cdot)}=\frac{1}{\lambda}p((1, I_{\Phi}(\lambda x))),\]
\[\parallel y\parallel_{\Phi,p(\cdot)}=\frac{1}{l}p((1, I_{\Phi}(ly))),\]
and define $q=\frac{2\lambda l}{\lambda+l}.$ Let us note that $\lambda x\neq ly$. Indeed, assuming that $\lambda x= ly$, we get by $\|x\|_{\Phi,p(\cdot)}=\|y\|_{\Phi,p(\cdot)}$ that $\lambda=l$ whence, by $x\neq y$, we obtain that $\lambda x\neq \lambda y$, a contradiction. Then by strict convexity of $I_\Phi$ and the fact that $p(\cdot)$ is strictly increasing on the half-line $\{(1,u):u\in R_+\}\subset R^2$, we get
\[\begin{aligned}
\parallel\frac{x+y}{2}\parallel_{\Phi, p(\cdot)}&\leq \frac{1}{q}p((1, I_{\Phi}(q\frac{x+y}{2})))\\
&=\frac{1}{q}p((1, I_{\Phi}(\frac{l}{l+\lambda}\lambda x+\frac{\lambda}{l+\lambda}ly)))\\
&<\frac{1}{q}p((1, \frac{l}{l+\lambda}I_{\Phi}(\lambda x)+\frac{\lambda}{l+\lambda}I_{\Phi}(ly)))\\
&=\frac{1}{q}p((\frac{l}{l+\lambda}+\frac{\lambda}{l+\lambda}, \frac{l}{l+\lambda}I_{\Phi}(\lambda x)+\frac{\lambda}{l+\lambda}I_{\Phi}(ly)))\\
&=\frac{1}{q}p(\frac{l}{l+\lambda}(1, I_{\Phi}(\lambda x)+\frac{\lambda}{l+\lambda}(1, I_{\Phi}(ly)))\\
&\leq\frac{1}{q}\{\frac{l}{l+\lambda}p((1,I_{\Phi}(\lambda x)))+\frac{\lambda}{l+\lambda}p((1, I_{\Phi}(ly)))\}\\
&=\frac{1}{2}\{\frac{1}{\lambda}p((1, I_{\Phi}(\lambda x)))+\frac{1}{l}p((1, I_{\Phi}(ly)))\}\\
&=\frac{1}{2}\{\parallel x\parallel_{\Phi,p(\cdot)}+\parallel y\parallel_{\Phi,p(\cdot)}\}\\
&=\frac{1}{2}\{1+1\}\\
&=1,
\end{aligned}\]
which finishes the proof that the positive cone $(L_{+}^{\Phi}(\mu),\parallel \cdot\parallel_{\Phi, p(\cdot)})$ is strictly convex. But then we obtain from a general result in [22] that the whole space  $(L_{\Phi}(\mu),\parallel \cdot\parallel_{\Phi, p(\cdot)})$ is strictly convex.\\
 \begin{COROLLARY} ([5]) Under the assumption that $\Phi$ is a strictly convex Orlicz function such that $\sup\limits_{u>0}[Au-\Phi(u)]=\infty,$ where $A:=\lim\limits_{u\rightarrow\infty}(\frac{\Phi(u)}{u})$, which gives that $K(x)\neq0$ for any $x\in L^{\Phi}(\mu)\setminus\{0\}$, the Orlicz space $L^{\Phi}(\mu)$ equipped with the Orlicz norm is strictly convex.
 \end{COROLLARY}
 Proof. Under the assumptions on $\Phi$, condition (2) from  Theorem 5 is satisfied (see [7]). Moreover, the Orlicz norm in $L^{\Phi}(\mu)$ is just the norm $\parallel \cdot\parallel_{\Phi, p(\cdot)}$ with $p((u,\nu))=\mid u \mid+\mid\nu\mid$ for any $(u, \nu)\in \mathbb{R}^{2}$. Since this norm $p(\cdot)$ is strictly increasing on the vertical half-line $\{(1, u): u\in \mathbb {R_{+}}\}$ in $R^{2}$, the thesis of our corollary follows directly form Theorem 5.\\
\begin{THEOREM} Assume that $p(\cdot)$ is a lattice norm in $\mathbb{R}^{2}$ which is strictly increasing on the vertical half-line $\{(1, u): u\in \mathbb {R_{+}}\}\subseteq \mathbb{R}^{2}$ and $\Phi$ is an Orlicz function such that $K(x)\neq \emptyset$ for any $x\in L^{\Phi}(\mu)\setminus\{0\}.$ Then the following statements are equivalent:\\
$(i)$   $a(\Phi)=0,$\\
$(ii)$  $(L^{\Phi}(\mu),\parallel \cdot\parallel_{\Phi, p(\cdot)})$ is strictly monotone,\\
$(iii)$ $(E^{\Phi}(\mu),\parallel \cdot\parallel_{\Phi, p(\cdot)})$ is strictly monotone.\\
\end{THEOREM}
{\it Proof.} $(i)\Longrightarrow(ii)$. Let $0\leq x\leq y\in S(L^{\Phi}(\mu),\parallel \cdot\parallel_{\Phi, p(\cdot)})$ and $x\neq y$. We know, by the assumption, that\\
\[\parallel y\parallel_{\Phi, p(\cdot)}=\frac{1}{l}p((1, I_{\Phi}(ly)))\]\\
for some $l\in(0,+\infty).$ Since $0\leq lx\leq ly$ and $lx\neq ly,$ and by the assumption that $a(\Phi)=0$ the Orlicz function $\Phi$ is strictly increasing on $\mathbb{R}_{+}$, we obtain\\
\[I_{\Phi}(lx)<I_{\Phi}(ly).\]\\
In consequence\\
\[\begin{aligned}
\parallel x\parallel_{\Phi, p(\cdot)}&\leq\frac{1}{l}p((1, I_{\Phi}(lx)))< \frac{1}{l}p((1, I_{\Phi}(ly)))\\
&=\parallel y\parallel_{\Phi, p(\cdot)},
\end{aligned}\]\\
which finishes the proof of the implication $(i)\Longrightarrow (ii)$. The implication $(ii)\Longrightarrow (iii)$ is obvious.\\
$(iii)\Longrightarrow (i)$. Assuming that $(i)$ does not hold, we will prove that $(iii)$ does not hold. Take any $y\in E^{\Phi}(\mu)$ such that $y\geq0$, $\parallel y\parallel_{\Phi, p(\cdot)}=1$ and the set $A:=\Omega\setminus suppy$ has positive measure. Define $z=y+\frac{a(\Phi)}{k}\mathcal{X}_{A},$ where $k\in K(y).$ Then $z\in E^{\Phi}(\mu)$ and $\parallel z\parallel_{\Phi, p(\cdot)}\geq1$ because $\mid z(t)\mid\geq\mid y(t)\mid$ for $\mu-a.e.$ $t\in \Omega.$ On the other hand $\parallel z\parallel_{\Phi, p(\cdot)}\leq\frac{1}{k}p((1, I_{\Phi}(kz)))=\frac{1}{k}p((1, I_{\Phi}(ky)))=\parallel y\parallel_{\Phi, p(\cdot)}=1.$ Therefore, $\parallel z\parallel_{\Phi, p(\cdot)}=1,$ whence $(iii)$ does not hold.\\
\begin{REMARK} Let us note that the assumption that the norm $p(\cdot)$ is strictly increasing on the half-line $\{(1, u): u\in \mathbb{R}_{+}\}$ is not necessary in general neither for strict convexity nor for strict monotonicity of the space $(L^{\Phi}(\mu),\parallel \cdot\parallel_{\Phi, p(\cdot)})$. Namely, if $p((u, \nu))=\max(\mid u\mid, \mid u\mid)({\forall}(u,\nu)\in \mathbb{R}^{2}),$ then $p(\cdot)$ is not strictly increasing on the half-line $\{(1,u):u\in \mathbb{R}_{+}\},$ but if $\Phi$ satisfies suitable $\Delta_{2}$-condition, then the space $(L^{\Phi}(\mu),\parallel \cdot\parallel_{\Phi, p(\cdot)})$ is strictly convex whenever $\Phi$ is strictly convex, and $(L^{\Phi}(\mu),\parallel \cdot\parallel_{\Phi, p(\cdot)})$ is strictly monotone, whenever $a(\Phi)=0$, see [33] and [25] (because $\parallel \cdot\parallel_{\Phi, p(\cdot)}$ is then equal to the Luxemburg norm).
\end{REMARK}
\begin{COROLLARY} It follows from Theorem 4 that if $\Phi$ is an Orlicz function with $a(\Phi)>0$ and if the measure space is non-atomic and infinite or the counting measure on $2^{N},$ then under the assumptions of Theorem 4 on the norm $p(\cdot)$, the space $(L^{\Phi}(\mu),\parallel \cdot\parallel_{\Phi, p(\cdot)})$ is not strictly monotone, because $(l^{\infty},\parallel\cdot\parallel_{\infty})$ is not strictly monotone (namely, the elements $x=(1,0,0,\cdots)$ and $y=(1,1,0,0,\cdots)$ satisfy the conditions $0\leq x\leq y,$ $x\neq y,$ $\parallel x\parallel_{\infty}=\parallel y\parallel_{\infty}=1$). Therefore, $(L^{\Phi}(\mu),\parallel \cdot\parallel_{\Phi, p(\cdot)})$ has no monotonicity property non-weaker than strict monotonicity. It is also not strictly convex (because any strictly convex Banach lattice is strictly monotone), so it has no convexity property non-weaker than strict convexity.
\end{COROLLARY}
\begin{THEOREM} Assume that $p(\cdot)$ is a lattice norm on $\mathbb{R}^{2}$ which is strictly monotone and $p((1,0))=p((0,1))=1$. Let $\Phi$ be an Orlicz function; $x,y\in L^{\Phi}(\mu),$ $0\leq x\leq y\in S((L^{\Phi}(\mu),\parallel \cdot\parallel_{\Phi, p(\cdot)}))$.\\
Then\\
(3) $\parallel y-x\parallel_{\Phi, p(\cdot)}\leq 1-\delta_{m, p(\cdot)}(I_{\Phi}(x)),$\\
where $\delta_{m, p(\cdot)}(\cdot)$ is the modulus of monotonicity of the space $(\mathbb{R}^{2},p(\cdot)).$\\
\end{THEOREM}
{\it Proof.} First let us note that since $\mathbb{R}^{2}$ is finitely dimensional, by strict monotonicity under the norm $p(.)$, $(\mathbb{R}^{2},p(\cdot))$ is uniformly monotone. Therefore, $\delta_{m, p(\cdot)}(\varepsilon)>0$ for any $\varepsilon\in(0,1)$. For the definition of the modulus of monotonicity and its properties see [23].\par
Let us take $x$ and $y$ mentioned in the theorem and any $\varepsilon>0.$ By the definition of the norm $\parallel \cdot\parallel_{\Phi, p(\cdot)},$ there exists $l\in (0,+\infty)$ such that\\
\[\frac{1}{l}p((1, I_{\Phi}(ly)))\leq\parallel y\parallel_{\Phi, p(\cdot)}+\varepsilon=1+\varepsilon.\]\\
Since $p((1, I_{\Phi}(ly)))\geq p((1, 0))=1,$ the previous inequality implies that $\frac{1}{l}\leq 1+\varepsilon,$ whence $l\geq\frac{1}{1+\varepsilon}.$ Hence, by the assumption that $p((0,1))=1,$ we have\\
(4) $p((0,\frac{1}{l}I_{\Phi}(lx)))=\frac{1}{l}I_{\Phi}(lx)p((0, 1))\geq(1+\varepsilon)I_{\Phi}(\frac{x}{1+\varepsilon}).$\\
By convexity of the modular $I_{\Phi(\cdot)}$, we obtain
\[0\leq(1+\varepsilon)I_{\Phi}(\frac{x}{1+\varepsilon})\leq{I_{\Phi}(x)}\leq{\parallel x\parallel_{\Phi, p(\cdot)}}\leq{\parallel y\parallel_{\Phi, p(\cdot)}}=1.\]\\
By superadditivity of $\Phi$ on $\mathbb{R}_{+}$, we have superadditivity of the modular $I_{\Phi}.$ Hence, and by the equality\\
\[p((0,\frac{1}{l}I_{\Phi}(lx)))=\frac{1}{l}I_{\Phi}(lx)p((0, 1))=\frac{1}{l}I_{\Phi}(lx),\]\\
we obtain\\
\[\begin{aligned}
\parallel y-x\parallel_{\Phi,p(\cdot)}&\leq\frac{1}{l}p((1,I_{\Phi}(l(y-x))))\\
&\leq\frac{1}{l}p((1,I_{\Phi}(ly)-I_{\Phi}(lx)))\\
&=p((\frac{1}{l},\frac{1}{l}I_{\Phi}(ly))-(0, \frac{1}{l}I_{\Phi}(lx)))\\
&\leq\frac{1}{l}p((1,I_{\Phi}(ly)))-\delta_{m,p(\cdot)}(p((0,\frac{1}{l}I_{\Phi}(lx))))\\
&\leq1+\varepsilon-\delta_{m,p(\cdot)}((1+\varepsilon)I_\Phi(\frac{x}{1+\varepsilon})).
\end{aligned},\]\\
Taking in place of $\varepsilon$ a sequence $\{\varepsilon_n\}_{n=1}^\infty$ such that $\varepsilon_n\searrow0$ as $n\nearrow\infty$ and applying the Beppo Levi theorem, we obtain that $(1+\varepsilon)I_\Phi(\frac{x}{1+\varepsilon})\nearrow I_\Phi(x)$ as $n\nearrow\infty$. Since the modulus of monotonicity is continuous on the interval $[0,1)$, we obtain the desired inequality.
\begin{REMARK} If $p(\cdot)$ is a lattice norm in $\mathbb{R}^{2}$ and $\Phi$ is an Orlicz function, then the Orlicz space $(L^{\Phi}(\mu),\|\cdot\|_{\Phi,p(\cdot)})$ is order continuous if and only if $\Phi$ satisfies suitable $\Delta_{2}$-condition .\\
Proof: Since all norms in $\mathbb{R}^{2}$ are equivalent, our norms $\|\cdot\|_{\Phi,p(\cdot)}$ are equivalent to the norm $\|\cdot\|_{\Phi,p(\cdot)}$ ,where $p_{0}((u,v))=\max (|u|,|v|).$ It is known (see [8] and [10]) that the norm $\|\cdot\|_{\Phi,p_{0}(\cdot)}$ is equal to the Luxemburg norm $\|\cdot\|_{\Phi}.$ It is also well known (see [5],[9]) that the Orlicz space $(L^{\Phi}(\mu),\|\cdot\|_{\Phi})$  is order continuous if and only if $\Phi$ satisfies suitable $\Delta_{2}$-condition. Since order continuity is preserved by equivalent norms, we obtain that the space $(L^{\Phi}(\mu),\|\cdot\|_{\Phi,p(\cdot)})$  is order continuous if and only if $\Phi$ satisfies suitable $\Delta_{2}$-condition.
\end{REMARK}
\begin{REMARK} It is known (see [5]) that for any sequence $\{x_{n}\}^{\infty}_{n=1}$ in $L^{\Phi}(\mu),$ with $I_{\Phi}(x_{n})\rightarrow 0$ as $n\rightarrow\infty$, we have $\|x_{n}\|_{\Phi}\rightarrow 0$ as $n\rightarrow\infty$ if and only if $a(\Phi)=0$ and $\Phi$ satisfies the suitable $\Delta_{2}$-condition. Since for any lattice norm $p(\cdot)$ in $\mathbb{R}^{2}$ with $p((1,0))=1$, the norm  $\|\cdot\|_{\Phi,p(\cdot)}$ is equivalent to the Luxemburg norm $\|\cdot\|_{\Phi}$ in $L^{\Phi}(\mu)$, we have the same dependence between $I_{\Phi}(x_{n})\rightarrow 0$ and $\|x_{n}\|_{\Phi,p(\cdot)}\rightarrow 0$ as  $n\rightarrow\infty$ if and only if $a(\Phi)=0$ and $\Phi$ satisfies suitable $\Delta_{2}$-condition. This is equivalent to the fact that
\begin{equation}\label{5}
(\forall \varepsilon >0)(\exists\delta(\varepsilon)>0)(\forall x\in{L^{\Phi}(\mu)} )(\|\cdot\|_{\Phi,p_{(\cdot)}}\geq\varepsilon\Rightarrow I_{\Phi}(x)\geq\delta(\varepsilon)).
\end{equation}
\end{REMARK}
\begin{THEOREM} Let $p(\cdot)$ be a strictly monotone lattice norm in $\mathbb{R}^{2}$ such that $p((1,0))= p((0,1))=1$ and $\Phi$ be an Orlicz function. Then the assumption $a(\Phi)=0$ implies that:\\
$(i)~(L^{\Phi}(\mu),\|\cdot\|_{\Phi,p(\cdot)})$ is strictly monotone, \\
$(ii)~(E^{\Phi}(\mu),\|\cdot\|_{\Phi,p(\cdot)})$ is lower locally uniformly monotone.\\
Assuming additionally that the couple $(\Phi,p(\cdot))$ is such that $K(x)\neq\emptyset$ for any $x\in E^{\Phi}(\mu)\backslash\{0\}$ the conditions $a(\Phi)=0$, $(i),(ii)$ and\\
$(iii)~(E^{\Phi}(\mu),\|\cdot\|_{\Phi,p(\cdot)})$ is strictly monotone,\\
are equivalent.\\
\end{THEOREM}
{\it Proof.} Under the assumption on $p(\cdot)$, we have by Theorem 7 that for any $x,y\in L^{\Phi}(\mu)$ such that $0\leq x\leq y,\|y\|_{\Phi,p(\cdot)}=1 $ and $\parallel x \parallel_{\Phi,p(\cdot)}\geq\varepsilon$, where $\varepsilon\in(0,1)$, we have

~~~~~~~~~~~~~~~~~~~~~~~~$\parallel y-x \parallel_{\Phi,p(\cdot)}\leq 1-\delta_{m,p(\cdot)} (I_{\Phi}(x)),$  \quad\quad\quad\quad\quad\quad\quad   (5)\\
where $\delta_{m,p(\cdot)}$ is the modulus of monotonicity of the space $(\mathbb{R}^{2},p(\cdot))$ which is uniformly monotone as a strictly monotone finite dimensional Banach space. Therefore, $\delta_{m,p(\cdot)}(I_{\Phi}(x))>0$ by the fact that $a(\Phi)=0$ implies that $I_{\Phi}(x)>0$. In consequence $\parallel y-x \parallel_{\Phi,p(\cdot)}< 1$, which means that property (i) holds.

Assume now that $x$ and $y$ are as above, but they belong to the space $E^{\Phi}(\mu)$. Let us note that the space $(E^{\Phi}(\mu), \parallel\parallel_{\Phi,p(\cdot)})$ has the following property.

$\forall y\in S_{+}((E^{\Phi}(\mu), \parallel \cdot \parallel_{\Phi,p(\cdot)})), \quad\forall \varepsilon\in(0,1), \quad\exists \delta(y,\varepsilon)>0,~such ~ that\\
 for ~ any ~  0\leq x\leq y:
 \parallel x \parallel_{\Phi,p(\cdot)}\geq\varepsilon\Rightarrow I_{\Phi}(x)\geq\delta(y,\varepsilon).~~~~~~~~~~~~~~~~~~~~~~~~~~~~~~~~~~~~~~~~~~~(6)$

Indeed, if this property does not hold, then there exist $y\in S_{+}((L^{\Phi}(\mu), \parallel \cdot \parallel_{\Phi,p(\cdot)}))$, $\varepsilon\in(0,1)$ and a sequence $\{x_{n}\}_{n=1}^{\infty}$ in  $E^{\Phi}(\mu)$ such that $\|x_n\|_{\Phi,p(\cdot)}\geq\varepsilon$, $0\leq x_{n}\leq y$ for any $n\in N$, and $I_{\Phi}(x_{n})\rightarrow 0$ as $n\rightarrow \infty$. However, the last condition implies that $x_{n}\rightarrow 0$ in measure as $n\rightarrow \infty$. Hence, by the assumption that the measure space $(\Omega, \Sigma, \mu)$ is $\delta$-finite, there exists a subsequence $\{x_{n_{k}}\}_{k=1}^{\infty}$ of $\{x_{n}\}_{n=1}^{\infty}$ such that $x_{n_{k}}\rightarrow 0$ as $k\rightarrow \infty$ $\mu-a.e.$ in $\Omega$. Hence $\Phi \circ \lambda x_{n_{k}} \rightarrow 0 $  $\mu-a.e.$ in $\Omega$ as $k\rightarrow\infty$ for any $\lambda >0$. Since, by the assumption that $y\in E^{\Phi}(\mu)$, we have $\Phi\circ\lambda y \in L^{1}(\mu)$ for any $\lambda >0$ and, by $0\leq x_{n}\leq y$ for any $n\in N$, we have $\Phi\circ\lambda x_{n_{k}}\leq \Phi\circ\lambda y$ for any $k\in N$ and $\lambda>0$, the Lebesgue dominated convergence theorem implies that $I_{\Phi}(\lambda x_{n_{k}})=\parallel\Phi \circ \lambda x_{n_{k}}\parallel_{L^{1}(\mu)}\rightarrow 0$ as $k\rightarrow\infty$ for any $\lambda>0$, which means that $\parallel x_{n_{k}} \parallel_{\Phi,p(\cdot)}\rightarrow 0$ as $k\rightarrow \infty$, a contradiction, which proves property (6) of $E^{\Phi}(\mu)$.
Conditions (5) and (6) yield
\begin{equation*}
\parallel y-x \parallel_{\Phi,p(\cdot)}\leq 1-\delta_{m,p(\cdot)}(\delta(y,\varepsilon)),
\end{equation*}
where $\delta(y,\varepsilon)$ does not depend on $x$, which means that property (ii) of $E^{\Phi}(\mu)$ holds.

It is obvious that (ii)$\Rightarrow$ (iii). In order to finish the proof, we need only to show that under the assumption that $K(x)\neq \emptyset$ for any $x\in E^{\Phi}(\mu)\backslash \{0\}$ condition (iii) implies that $a(\Phi)=0$. Assume that $a(\Phi)>0$ and $K(x)\neq \emptyset$ for any $x\in E^{\Phi}(\mu)\backslash \{0\}$.
Let us take any $x\in E^{\Phi}(\mu)$ such that $x\geq 0$, $\parallel x \parallel_{\Phi,p(\cdot)}=1$ and $\mu(A)>0$, where $A=\Omega\backslash suppx$.
Defining $y=x+\dfrac{a(\Phi)}{k} \mathcal{X}_{A}$, where $k\in K(x)$, we have $0\leq x\leq y\in E^{\Phi}(\mu)$ and $x\neq y$.
Hence $\parallel y\parallel_{\Phi,p(\cdot)}\geq 1$. On the other hand, by $I_{\Phi}(ky)=I_{\Phi}(kx)$, we have
\begin{equation*}
\parallel y\parallel_{\Phi,p(\cdot)}\leq \dfrac{1}{k}p((1,I_{\Phi}(ky)))=\dfrac{1}{k}p((1,I_{\Phi}(kx)))=\parallel x\parallel_{\Phi,p(\cdot)}=1.
\end{equation*}
Therefore, property (iii) of $(E^{\Phi}(\mu),\parallel \cdot \parallel_{\Phi,p(\cdot)})$ does not hold if $a(\Phi)=0$, and the proof of our theorem is complete.

\begin{THEOREM}
Let $p(\cdot)$ be a strictly monotone lattice norm on $\mathbb{R}^{2}$ satisfying the condition $p((1,0))=p((0,1))=1$ and $\Phi$ be an Orlicz function. Consider the following conditions:\\
(i) $a(\Phi)=0$ and $\Phi$ satisfied suitable $\Delta_{2}$-condition,\\
(ii) $(L^{\Phi}(\mu),\parallel\cdot\parallel_{\Phi,p(\cdot)})$ is uniformly monotone,\\
(iii) $(L^{\Phi}(\mu),\parallel\cdot\parallel_{\Phi,p(\cdot)})$ is upper locally uniformly monotone,\\
(iv) $(E^{\Phi}(\mu),\parallel\cdot\parallel_{\Phi,p(\cdot)})$ is upper locally uniformly monotone.\\
Then (i)$\Rightarrow$ (ii) $\Rightarrow$ (iii) $\Rightarrow$ (iv). Assuming additionally that $K(x)\neq \emptyset$ for any $x\in E^{\Phi}(\mu)\backslash\{0\}$,
we have also that (iv)$\Rightarrow$ (i), whence we have then that all these four conditions are equivalent.
\end{THEOREM}
{\it Proof.}
(i) $\Rightarrow$ (ii). Assume that condition (i) is satisfied and that $0\leq x\leq y\in S((L^{\Phi}(\mu),\parallel\cdot\parallel_{\Phi,p(\cdot)}))$ and
$\parallel x\parallel_{\Phi,p(\cdot)}\geq\varepsilon$. By the assumption that $\Phi$ satisfies suitable $\Delta_{2}-condition$, we know that for any sequence $\{x_{n}\}_{n=1}^\infty$ in $L^{\Phi}(\mu)$, the conditions $I_{\Phi}(x_{n})\rightarrow 0$ as $n\rightarrow \infty$ and $\parallel x_{n}\parallel_{\Phi,p(\cdot)}\rightarrow 0$ as $n\rightarrow \infty$ are equivalent.~Therefore, there exists $\delta(\varepsilon)\in(0,1)$ such that $I_{\Phi}(x)\geq\delta(\varepsilon)$. By inequality (3) from Theorem $7$, we have
\begin{equation*}
\parallel y-x\parallel_{\Phi,p(\cdot)}\leq 1-\delta_{m,p(\cdot)}(I_{\Phi}(x))\leq 1-\delta_{m,p(\cdot)}(\delta(\varepsilon)).
\end{equation*}
Since $\delta_{m,p(\cdot)}(\delta(\varepsilon))\in(0,1)$, property (ii) holds.

The implications (ii) $\Rightarrow$ (iii) $\Rightarrow$ (iv) are obvious, so in order to finish our proof, we need only to prove that (iv) $\Rightarrow$ (i). Assume that condition (i) is not satisfied and $K(x)\neq \emptyset$ for any $x\in E^{\Phi}(\mu)\backslash \{0\}$.
Then we have the alternative $a(\Phi)=0$ and $K(x)\neq0$ for any $x\in E^{\Phi}(\mu)\backslash\{0\}$ or
$\Phi\notin\Delta_{2}$ and $K(x)\neq \emptyset$ for any $x\in E^{\Phi}(\mu)\backslash\{0\}$.

In the first situation, by virtue of Theorem $5$, condition (iv) does not hold, because $(E^{\Phi}(\mu),\parallel\cdot\parallel_{\Phi,p(\cdot)})$ is then not strictly monotone.

In the second situation take a set $A\in\Sigma$ such that $0<\mu(A)<\mu(\Omega)$. By Proposition $2.1$ in [8] there exists a sequence
$\{y_{n}\}_{n=1}^\infty$ in $E^{\Phi}(\mu)$ such that $supp y_{n}\subset T\backslash A$,
$\dfrac{1}{1+2^{-n}}\leq\parallel y_{n}\parallel_{\Phi}$ and $I_{\Phi}(y_{n})\leq 2^{-n}$ for any $n\in N$.
Then $\parallel y_{n}\parallel_{\Phi,p(\cdot)}\geq\parallel y_{n}\parallel_{\Phi}\geq\dfrac{1}{1+2^{-n}}$ for any $n\in N$.
Take any $x\in E^{\Phi}(\mu)$ such that $suppx\subset A$ and $\parallel x\parallel_{\Phi,p(\cdot)}=1$.
Let $k\in K(x)$. Then $k>1$, because for any $k\in(0,1]$ we have
\begin{equation*}
\dfrac{1}{k}p((1,I_{\Phi}(kx)))>\dfrac{1}{k}p((1,0))=\dfrac{1}{k}=1=\parallel x\parallel_{\Phi,p(\cdot)},
\end{equation*}
which means that $k\notin K(x)$. Defining $x_{n}=(\dfrac{1}{k})y_{n}$, we have
$\|x_n\|_{\Phi,p(\cdot)}=(\dfrac{1}{k})\|y_{n}\|\geq\dfrac{1}{k(1+2^{-n})}\geq\dfrac{2}{3k}$
and $I_{\Phi}(kx_{n})=I_{\Phi}(y_{n})\leq 2^{-n}$. Hence
\begin{equation*}
\begin{split}
\parallel x+x_{n}\parallel_{\Phi,p(\cdot)}&\leq \dfrac{1}{k}p((1,I_{\Phi}(k(x+x_{n}))))=\dfrac{1}{k}p((1,I_{\Phi}(kx))+I_{\Phi}(kx_{n})))\\
&=\dfrac{1}{k}p((1,I_{\Phi}(kx))+(0,I_{\Phi}(kx_{n})))\\
&\leq\dfrac{1}{k}p((1,I_{\Phi}(k(x)))+\dfrac{1}{k}p((0,I_{\Phi}(k(x_{n})))\\
&\leq\parallel x\parallel_{\Phi,p(\cdot)}+\dfrac{1}{k}p((0,2^{-n}))\\
&\leq 1+2^{-n},
\end{split}
\end{equation*}
which together  with $\|x_n\|_{\Phi,p(\cdot)}\geq\dfrac{2}{3k}$ for any $n\in N$ means that the space $(E^{\Phi}(\mu),
\parallel\cdot\parallel_{\Phi,p(\cdot)})$
is not locally upper uniformly monotone. This finishes the proof.




\end{document}